\documentclass[12pt]{article}
\usepackage{a4wide}
\usepackage{theorem}
\usepackage{amsfonts}
\usepackage{amsmath}
\usepackage{times}

\theorembodyfont{\sl}
\newtheorem{thm}[section]{Theorem}
\newtheorem{prop}[section]{Proposition}

\theorembodyfont{\rmfamily}
\newtheorem{defi}[section]{Definition}

\newenvironment{prf}[1]{\trivlist
\item[\hskip \labelsep{\it
#1.\hspace*{.3em}}]}{~\hspace{\fill}~$\square$\endtrivlist}

\newenvironment{proof}{\begin{prf}{\bf Proof}}{\end{prf}}

\newcommand{\ol}{\overline}
\newcommand{\ZZ}{{\mathbb Z}}
\newcommand{\QQ}{{\mathbb Q}}
\newcommand{\CC}{{\mathbb C}}
\newcommand{\FF}{{\mathbb F}}
\newcommand{\Qbar}{{\overline{\QQ}}}
\newcommand{\Fbar}{{\overline{\FF}}}
\DeclareMathOperator{\Spec}{Spec}
\newcommand{\calE}{{\cal E}}
\newcommand{\End}{{\rm End}}
\newcommand{\Aut}{{\rm Aut}}
\newcommand{\et}{\mathrm{et}}

\begin{document}
\title{On the $p$-adic geometry of traces of singular moduli} 
\author{Bas Edixhoven}
\maketitle

The aim of this article is to show that $p$-adic geometry of modular
curves is useful in the study of $p$-adic properties of \emph{traces}
of singular moduli. In order to do so, we partly answer a question by
Ono (\cite[Problem~7.30]{Ono1}). As our goal is just to illustrate how
$p$-adic geometry can be used in this context, we focus on a
relatively simple case, in the hope that others will try to obtain the
strongest and most general results. For example, for $p=2$, a result
stronger than Thm.~\ref{thm1} is proved in~\cite{Boylan1}, and a
result on some modular curves of genus zero can be found
in~\cite{Osburn1} . It should be easy to apply our method, because of
its local nature, to modular curves of arbitrary level, as well as to
Shimura curves.

\begin{defi}\label{defi1}
For $d$ a positive integer that is congruent to $0$ or~$3$ mod~$4$,
let $O_d$ be the quadratic order of discriminant~$-d$. For $f$
in~$\ZZ[j]$ and $E$ an elliptic curve over some ring~$R$, let $f(E)$
be the element of $R$ obtained by evaluating the $f$ on the
$j$-invariant of~$E$. For such $d$ and~$f$, let:
$$
t_f(d) := \sum_{\End(E)\supset O_d} 2f(E)/\#\Aut(E),
$$
where the sum ranges over the set of isomorphism classes of complex
elliptic curves whose ring of endomorphisms contains~$O_d$. We also
define an integer $\alpha(d)$ to be~$2$ if
$\QQ(\sqrt{-d})=\QQ(\sqrt{-1})$, $3$ if
$\QQ(\sqrt{-d})=\QQ(\sqrt{-3})$, and~$1$ otherwise.
\end{defi}
For $d$ as above and $m$ a positive integer, the number $t_m(d)$
defined in~\cite{Ono1} is obtained by taking $f:=(j-744)|T_0(m)$. 

\begin{thm}\label{thm1}
Let $d>0$ be an integer that is congruent to $0$ or~$3$ modulo~$4$,
and let $f$ be in~$\ZZ[j]$. Let $p$ be a prime not dividing~$d$ that
splits in~$\QQ(\sqrt{-d})$, and let $n\geq1$. Then
$\alpha(d)t_f(p^{2n}d)$ is an integer, and
$\alpha(d)t_f(p^{2n}d)\equiv 0$ mod~$p^n$.
\end{thm}
All that we need from the local moduli theory of ordinary elliptic
curves in positive characteristic is summarized in the following
proposition. Definitions for the terms occurring in it can be found in
\cite{Silverman1} and~\cite{Katz1}.
\begin{prop}\label{prop2}
Let $p$ be a prime number, $k$ a finite field of
characteristic~$p$. Let $E_0$ be an ordinary elliptic curve over~$k$,
and let $A$ be its endomorphism ring. Then $A$ is an order in a
quadratic extension of~$\QQ$, and $A$ is split at~$p$: $\ZZ_p\otimes
A$ is isomorphic to $\ZZ_p\times \ZZ_p$; in particular, $A$ is maximal
at~$p$.

Let $k\to \ol{k}$ be an algebraic closure, and let $W$ be the ring of
Witt vectors of~$k$. Let $E/R$ be the universal deformation
of~$E_0/\ol{k}$ over $W$-algebras. Let $\alpha\colon \QQ_p/\ZZ_p\to
E_0(\ol{k})[p^\infty]$ be a trivialisation of the group of torsion
points of $p$-power order of~$E_0(\ol{k})$. Then $\alpha$ induces a
so-called Serre-Tate parameter $q\in R^*$, and $R=W[[q-1]]$.  For
$n\geq 0$, let $A_n$ be the subring $\ZZ+p^nA$ of~$A$, i.e., the order
of index $p^n$ in~$A$. Then the closed subscheme of $\Spec R$ over
which all elements of~$A_n$ lift as endomorphisms of~$E$ is the closed
subscheme defined by the equation $q^{p^n}=1$.
\end{prop}
\begin{proof}
The endomorphism ring~$A$ is free of finite rank as a $\ZZ$-module,
and is an integral domain because each non-zero element in it is
surjective as a morphism from $E_0$ to itself. The $p$-divisible group
$E_0[p^\infty]$ is the direct sum of its local and etale parts
$E_0[p^\infty]^0$ and $E_0[p^\infty]^\et$, hence its endomorphism ring
is the $\ZZ_p$-algebra $\ZZ_p\times\ZZ_p$. Then $A$ is commutative
because it embeds into $\ZZ_p\times\ZZ_p$. The image of the Frobenius
endomorphism of~$E_0$ is of the form $(p^mu,v)$, with $u$ and~$v$
in~$\ZZ_p^*$, and $|k|=p^m$. This proves that $A$ is quadratic
over~$\ZZ$, and split at~$p$.

The construction of~$q$ and the statement that $R=W[[q-1]]$ are
in~\cite[\S2]{Katz1}. There it is also shown every $f$ in~$A$
determines a closed subscheme $V_f$ of~$\Spec R$ given by the
condition that $f$ can be lifted as an endomorphism of~$E_{V_f}$, and
universal for that property. This subscheme $V_f$ coincides with the
closed subscheme over which $f$ can be lifted as an endomorphism
of~$E[p^\infty]$. Let $n\geq0$, and let $V$ be the intersection of the
$V_f$ for all $f$ in~$A_n$. Then $V$ is defined by the condition that
the endomorphism $(p^n,0)$ of~$E_0[p^\infty]$ lifts as an endomorphism
of the $p$-divisible group. The proof of Part~4 of
\cite[Thm~2.1]{Katz1} shows that $V$ is defined by the
equation~$q^{p^n}=1$.
\end{proof} 

We can now prove Theorem~\ref{thm1}. Let $d$, $f$, $p$ and~$n$ be as
in the statement. For each $E$ with $\End(E)\supset O_{p^{2n}d}$ we
have that $\#\Aut(E)/2$ divides~$\alpha(d)$, hence
$\alpha(d)t_f(p^{2n}d)$ is an integer. Let $\Qbar$ be the algebraic
closure of~$\QQ$ in~$\CC$, and choose an embedding of~$\Qbar$ into an
algebraic closure~$\Qbar_p$ of~$\QQ_p$. For $E$ an elliptic curve
over~$\CC$ with complex multiplications, let $\ol{E}$ denote its
reduction over~$\Fbar_p$. For each $E$ in the sum in
Definition~\ref{defi1} there is a unique $E_0$ such that
$\End(E_0)=O_d$ and $\ol{E}\cong\ol{E_0}$. For each such~$E_0$, let
$\calE(E_0)$ denote the set of~$E$ with $\End(E)\supset O_d$
and~$\ol{E}\cong\ol{E_0}$. Then we have:
$$
\alpha(d)t_f(p^{2n}d) = 
\sum_{\End(E_0)=O_d}\frac{2\alpha(d)}{\#\Aut(E_0)}\#\Aut(E_0)
\sum_{E\in\calE(E_0)} f(E)/\#\Aut(E). 
$$
By construction, the $2\alpha(d)/\#\Aut(E_0)$ are integers. Fix an
$E_0$ as in the sum, and let $q$ be a $q$-parameter of the deformation
space of~$\ol{E_0}$. Then the relation between deformation spaces and
coarse moduli spaces (see~\cite[I~\S8.2.1]{DeRa}) and
Proposition~\ref{prop2} imply that:
$$
\#\Aut(E_0)\sum_{E\in\calE(E_0)} f(E)/\#\Aut(E) = 
\sum_{x^{p^n}=1}f(x-1),
$$ 
where we can now view $f$ as an element $\sum_{k\geq0}f_kt^k$
of~$W[[t]]$, with $t=q-1$ and $W$ the ring of Witt vectors
of~$\Fbar_p$. The observation that $\sum_{x^{p^n}=1}(x-1)^k$ is
in~$p^nW$ for all $k\geq0$ finishes the proof.

\vfill
\noindent
Bas Edixhoven\\
Universiteit Leiden\\
Mathematisch Instituut\\
Postbus 9512\\
2300\ RA\ \ Leiden\\
The Netherlands

\medskip\noindent
edix@math.leidenuniv.nl

\end{document}